\documentclass[11pt,a4paper]{article}

\usepackage[T1]{fontenc}
\usepackage[utf8]{inputenc}
\usepackage[english]{babel}
\usepackage{lmodern}
\usepackage{microtype}
\usepackage{geometry}
\usepackage{amsmath,amssymb,amsthm,mathtools,bbm}
\usepackage{graphicx}
\usepackage{float}
\usepackage[hidelinks]{hyperref}
\usepackage[nameinlink,capitalise,noabbrev]{cleveref}
\crefname{assumption}{assumption}{assumptions}
\Crefname{assumption}{Assumption}{Assumptions}

\theoremstyle{plain}
\newtheorem{theorem}{Theorem}[section]
\newtheorem{lemma}[theorem]{Lemma}
\newtheorem{proposition}[theorem]{Proposition}
\newtheorem{corollary}[theorem]{Corollary}

\theoremstyle{definition}

\newtheorem{assumption}[theorem]{Assumption}

\theoremstyle{remark}
\newtheorem{remark}[theorem]{Remark}

\numberwithin{equation}{section}

\newcommand{\R}{\mathbb{R}}
\newcommand{\E}{\mathbb{E}}
\newcommand{\Var}{\mathrm{Var}}
\newcommand{\1}{\mathbbm{1}}
\newcommand{\h}{h}

\title{\textbf{The Weighted Kolmogorov Metric Does Not Restore the Berry--Esseen Rate:} \\
\large Exact Rates for Heavy-Tailed Sums and an Impossibility Theorem}
\author{Armen Petrosyan}
\date{\today}

\begin{document}
\maketitle

\begin{abstract}
An earlier version of this note (January 2026) claimed that weighting the Kolmogorov
distance by $w_q(t)=(1+\h(t))^{-q}$ restores the Berry--Esseen rate $O(n^{-1/2})$ for
i.i.d.\ sums under the sub-cubic moment condition $\E|X|^{2+\delta}<\infty$. That claim
is false and is withdrawn, for two independent reasons.

First, the parameter region of the truncation scheme underlying the claimed proof is
empty: making all three terms of the core/tail decomposition $O(n^{-1/2})$ forces
$\delta=1$, i.e.\ a finite third moment, in which case the classical unweighted bound
already applies.

Second, and irreparably, we prove a lower bound. For symmetric returns with regularly
varying tails of index $\alpha\in(2,3)$, \emph{every} metric of the form
$\sup_t w(t)\,|F(t)-G(t)|$ with $1/w$ of polynomial growth satisfies
\[
\sup_t w(t)\,\big|P(Z_n\le t)-\Phi(t)\big|\ \ge\ c\, n^{-(\alpha-2)/2},
\qquad Z_n=\tfrac{1}{\sigma\sqrt n}\textstyle\sum_{i\le n}(X_i-\mu).
\]
Combined with the classical upper bound under $2+\delta$ moments, this identifies the
exact rate $n^{-(\alpha-2)/2+o(1)}$ of the weighted metric --- the same as the uniform
Kolmogorov distance. Reweighting in the argument $t$ is \emph{rate-neutral}: the
Berry--Esseen defect under heavy tails is located at bounded argument, where any
reasonable weight is of order one.

We record what does restore the $n^{-1/2}$ rate (correcting the Gaussian target by the
stable correction term; trimming the data rather than the metric), which parts of the
original proposal survive (the weighted statistic as a bootstrap-calibrated diagnostic),
and why the numerical experiments of the January version were dominated by Monte Carlo
noise. An appendix lists, statement by statement, the corrections with respect to that
version.
\end{abstract}

\section{Introduction}

\subsection{The claim, and the verdict}

Let $X_1,\dots,X_n$ be i.i.d.\ with $\mu=\E X_1$, $\sigma^2=\Var(X_1)\in(0,\infty)$, and
$Z_n=\frac{1}{\sigma\sqrt n}\sum_{i=1}^n(X_i-\mu)$. When $\E|X_1|^3<\infty$, the
Berry--Esseen theorem gives
$\sup_t|P(Z_n\le t)-\Phi(t)|=O(n^{-1/2})$. When only
$\E|X_1|^{2+\delta}<\infty$ with $\delta\in(0,1)$ --- the regime relevant to
heavy-tailed financial returns --- the uniform rate degrades to $O(n^{-\delta/2})$
\cite{Bikelis1966,Petrov1995}, and this is sharp in general.

The January 2026 version of this note proposed to repair the rate by changing the
\emph{metric}: for an exhaustion function $\h:\R\to[0,\infty)$ with $\h(t)\asymp|t|$ at
infinity and $q>0$, define
\begin{equation}\label{eq:metric}
d_{K,\h,q}(F,G):=\sup_{t\in\R}\,(1+\h(t))^{-q}\,|F(t)-G(t)|,
\end{equation}
and it claimed that, for regularly varying tails of index $\alpha>2+\delta$, a suitable
choice of $q$ yields $d_{K,\h,q}(\mathcal L(Z_n),\Phi)=O(n^{-1/2})$.

This is false. The purpose of the present version is to withdraw the claim and to
replace it by the correct picture, which consists of three statements.

\begin{itemize}
\item[(i)] \emph{The proof could not have worked} (\cref{sec:empty}): in the core/tail
truncation scheme, the Berry--Esseen term on the truncated core is of order
$(1+R_n)^{1-\delta}/\sqrt n$, so requiring it to be $O(n^{-1/2})$ forces the truncation
level $R_n$ to be bounded --- incompatible with the two remaining terms unless
$\delta=1$. The admissible parameter region is empty whenever $\E|X|^3=\infty$.
\item[(ii)] \emph{No proof can work} (\cref{sec:lower}): we prove that for symmetric
regularly varying tails of index $\alpha\in(2,3)$, every weight $w$ with $1/w$ of
polynomial growth satisfies
$\sup_t w(t)|P(Z_n\le t)-\Phi(t)|\ge c\,n^{-(\alpha-2)/2}$. The leading CLT correction
under heavy tails is a fixed profile of size $n^{1-\alpha/2}$ living at \emph{bounded}
argument, where the weight is inactive. Since the classical bound gives the matching
upper rate, the exact rate of \eqref{eq:metric} equals that of the uniform metric:
weighting in $t$ changes constants, never the rate.
\item[(iii)] \emph{What actually restores $n^{-1/2}$} (\cref{sec:whatworks}): one must
modify either the \emph{target} (compare $\mathcal L(Z_n)$ to the Gaussian plus its
stable correction term, an Edgeworth-type expansion in the heavy-tailed regime
\cite{IbragimovLinnik,ChristophWolf}) or the \emph{statistic} (trimmed or winsorized
sums \cite{CsorgoHaeuslerMason}). Ironically, the January version dismissed
winsorization in favour of the smooth metric weight; the correct moral is the opposite:
truncating the \emph{data} is exactly the operation that repairs the rate, and it is
what the (valid part of the) truncation proof was secretly doing, while the metric-side
weight contributed nothing.
\end{itemize}

What survives of the original proposal is operational, not asymptotic: the weighted
statistic \eqref{eq:metric} remains a legitimate goodness-of-fit diagnostic with
parametric-bootstrap critical values, whose tail-insensitivity is a finite-sample
variance statement rather than a rate statement. \Cref{sec:practice} records this in
corrected form. \Cref{sec:numerics} explains why the simulations of the January version
(160 Monte Carlo batches) could not have detected any of the rates discussed, and gives
a deterministic protocol. \Cref{app:corrections} lists the retractions individually.

\subsection{Setting and standing assumptions}

Throughout, $\Phi$ is the standard normal cdf, $\1\{\cdot\}$ denotes indicators, and
$F_n(t):=P(Z_n\le t)$.

\begin{assumption}[Moments]\label{ass:moment}
$\E|X_1-\mu|^{2+\delta}<\infty$ for some $\delta\in(0,1]$, and
$\sigma^2\in(0,\infty)$.
\end{assumption}

\begin{assumption}[Regular exhaustion]\label{ass:h}
$\h:\R\to[0,\infty)$ is Borel, locally bounded, $\h(t)\to\infty$ as $|t|\to\infty$, and
there are $c_1,c_2>0$, $t_0\ge0$ with $c_1|t|\le \h(t)\le c_2|t|$ for $|t|\ge t_0$.
\end{assumption}

(The local boundedness, implicit in the January version's use of
$\max_{|t|\le R}\h(t)$, is now explicit; it is needed both there and in
\cref{thm:lower}.) The following two facts from the January version are correct and are
kept without change.

\begin{proposition}[Metric property]\label{prop:metric}
For any $q>0$ and any finite $\h$, $d_{K,\h,q}$ is a metric on the set of cdfs.
\end{proposition}

\begin{lemma}[Weight equivalence]\label{lem:weight-equiv}
Under \cref{ass:h}, for every $q>0$ there are $C_-,C_+>0$ with
$C_-(1+|t|)^{-q}\le(1+\h(t))^{-q}\le C_+(1+|t|)^{-q}$ on $\R$; hence
$d_{K,\h,q}$ and $d_{K,|\cdot|,q}$ are equivalent, and so are the metrics attached to
two coarsely equivalent exhaustions.
\end{lemma}

By \cref{lem:weight-equiv} we may and do work with $w_q(t)=(1+|t|)^{-q}$ in proofs; all
statements transfer to general $\h$ with adjusted constants.

\section{The truncation route cannot work: the parameter region is empty}
\label{sec:empty}

The January version's argument ran as follows. Truncate at level $R$, apply a
(non-uniform) Berry--Esseen bound to the centered truncated variables, control the
difference between the full and truncated sums by the tail moment, and pay
$(1+R)^{-q}$ for the region $\{\h>R\}$ where only the weight is used. Writing
$M_3(R):=\E[|X-\mu|^3\1\{\h(X)\le R\}]\le C_\delta(1+R)^{1-\delta}\E|X-\mu|^{2+\delta}$
(interpolation on the core) and $\tau_R^2$ for the truncated variance
($\tau_R\to\sigma$), the resulting bound with $R_n=n^\beta$ has the three exponents
\begin{equation}\label{eq:threeterms}
n^{\beta(1-\delta)-\frac12},\qquad n^{-\beta\eta},\qquad n^{-\beta q},
\end{equation}
where $\eta>0$ governs the tail remainder
$\E[|X-\mu|^{2+\delta}\1\{\h(X)>R\}]\le KR^{-\eta}$.

\begin{proposition}[Empty admissible region]\label{prop:empty}
All three terms in \eqref{eq:threeterms} are $O(n^{-1/2})$ if and only if
\[
\beta(1-\delta)\le 0,\qquad \beta\eta\ge\tfrac12,\qquad \beta q\ge\tfrac12 .
\]
The last two force $\beta>0$; the first then forces $\delta=1$, i.e.\
$\E|X-\mu|^3<\infty$. In particular, whenever $\E|X-\mu|^3=\infty$ --- the entire
regime the construction was designed for --- no choice of $(\beta,q,\delta,\eta)$
yields the rate $n^{-1/2}$ through this scheme.
\end{proposition}

\begin{proof}
Immediate from \eqref{eq:threeterms}: $n^{\beta(1-\delta)-1/2}=O(n^{-1/2})$ iff
$\beta(1-\delta)\le0$.
\end{proof}

\begin{remark}[The error in the January version]\label{rem:arith}
Proposition~4.2 there imposed $\beta(1-\delta)\le\frac12$, which makes the first term
$O(1)$, not $O(n^{-1/2})$; Theorem~3.9 imposed no constraint on the first term at all.
The advertised ``balanced choice'' $\beta^\star=1/(2\eta)$, $q^\star=\eta$ is
disastrous when $\eta$ is small: for Student-$t$ with $\nu=2.5$ and $\delta=0.4$
(so $\eta=\nu-2-\delta=0.1$), it gives $\beta^\star=5$ and a ``core'' term
$n^{5(1-0.4)-1/2}=n^{2.5}\to\infty$ --- the bound is vacuous precisely on the paper's
own headline example.
\end{remark}

\begin{remark}[A second gap in the decomposition lemma]\label{rem:missing-n}
Independently of the exponent arithmetic, the truncation remainder in Lemma~4.1 of the
January version omitted a factor $n$: replacing $S_n$ by the truncated sum costs at
least $P(\exists i\le n:\h(X_i)>R)\approx n\,P(\h(X)>R)$, plus a mean-shift term of
order $\sqrt n\,|\E[(X-\mu)\1\{\h(X)>R\}]|$, neither of which appears in the stated
bound with ``absolute constants''. With $R_n=n^\beta$ these extra terms are
$n^{1-\beta\alpha}$ and $n^{\frac12-\beta(\alpha-1)}$ under a tail of index $\alpha$;
they are harmless for large $\beta$ but must be present. We do not repair the lemma in
detail, since by \cref{prop:empty} the route is dead regardless, and since the
\emph{correct} upper bound is both stronger and shorter (\cref{thm:upper}).
\end{remark}

\section{The exact rate of the weighted metric}\label{sec:lower}

\subsection{Upper bound: one line from the classical theory}

\begin{theorem}[Upper bound]\label{thm:upper}
Under \cref{ass:moment,ass:h}, for every $q\ge0$,
\[
d_{K,\h,q}\big(\mathcal L(Z_n),\Phi\big)\ \le\ \sup_{t\in\R}|F_n(t)-\Phi(t)|
\ \le\ \frac{C\,\E|X-\mu|^{2+\delta}}{\sigma^{2+\delta}\,n^{\delta/2}} ,
\]
with $C=C(\delta)$ absolute. Moreover, by the non-uniform bound of Bikelis
\cite{Bikelis1966} (see also \cite[Ch.~V]{Petrov1995}),
$|F_n(t)-\Phi(t)|\le C\,\E|X-\mu|^{2+\delta}\,\sigma^{-(2+\delta)}
n^{-\delta/2}(1+|t|)^{-(2+\delta)}$, so the weight improves the \emph{constant}
attached to the tail region --- and nothing else.
\end{theorem}

\begin{proof}
The first inequality is $w_q\le1$; the second is the Katz--Petrov extension of
Berry--Esseen under $2+\delta$ moments \cite{Bikelis1966,Petrov1995}.
\end{proof}

\subsection{Lower bound: the defect lives at bounded argument}

We now show the rate $n^{-\delta/2}$ (more precisely $n^{-(\alpha-2)/2}$ under a tail
index $\alpha$) cannot be improved by \emph{any} reweighting in $t$. We work in the
symmetric case, which covers Student-$t$ and symmetrized Pareto; the asymmetric case is
analogous with a complex tail constant, see \cref{rem:asym}.

\begin{assumption}[Symmetric regularly varying tails]\label{ass:RV}
$X$ is symmetric, $\sigma^2\in(0,\infty)$, and there are $\alpha\in(2,3)$, $c_0>0$ with
\[
P(|X|>x)=c_0\,x^{-\alpha}\big(1+\varepsilon(x)\big),\qquad
\varepsilon(x)\to0 \text{ as } x\to\infty .
\]
\end{assumption}

Student-$t$ with $\nu\in(2,3)$ satisfies \cref{ass:RV} with $\alpha=\nu$. Note that
\cref{ass:RV} implies $\E|X|^{2+\delta}<\infty$ exactly for $\delta<\alpha-2$.

\begin{lemma}[Characteristic function expansion]\label{lem:cf}
Under \cref{ass:RV}, as $t\to0$,
\[
1-\varphi(t)=\frac{\sigma^2t^2}{2}-c_0\,\kappa_\alpha\,|t|^\alpha+o(|t|^\alpha),
\qquad
\kappa_\alpha:=\int_0^\infty \frac{u-\sin u}{u^{\alpha}}\,du\ \in(0,\infty).
\]
\end{lemma}

\begin{proof}
Since $g(u):=\cos u-1+u^2/2\ge0$ with $g(u)\sim u^4/24$ at $0$ and $g(u)\sim u^2/2$ at
infinity, and since $X$ is symmetric,
\[
1-\varphi(t)=\int(1-\cos tx)\,dF(x)
=\frac{\sigma^2t^2}{2}-J(t),\qquad
J(t):=2\int_0^\infty g(tx)\,dF(x)\ \ge0 .
\]
Let $H(x):=P(X>x)$. Integrating by parts (the boundary terms vanish: $g(tx)H(x)\to0$ at
$0$ because $g(u)=O(u^4)$, and at $\infty$ because $g(tx)H(x)=O(x^{2-\alpha})$ for
$\alpha>2$),
\[
J(t)=2t\int_0^\infty g'(tx)\,H(x)\,dx,\qquad g'(u)=u-\sin u\ \ge 0,\ \text{non-decreasing}.
\]
Fix $x_1$ with $|\varepsilon|\le M<\infty$ on $[x_1,\infty)$. The compact part is
$2t\int_0^{x_1}g'(tx)H(x)\,dx\le 2t\,g'(tx_1)\,x_1=O(t^4)=o(|t|^\alpha)$. On the tail,
$H(x)=\tfrac{c_0}{2}x^{-\alpha}(1+\varepsilon(x))$, and the substitution $u=tx$ gives
\[
2t\int_{x_1}^\infty g'(tx)H(x)\,dx
=c_0\,|t|^\alpha\int_{tx_1}^\infty \frac{g'(u)}{u^{\alpha}}
\big(1+\varepsilon(u/|t|)\big)\,du .
\]
As $t\to0$, $\varepsilon(u/|t|)\to0$ pointwise for $u>0$, and the integrand is
dominated by $(1+M)\,g'(u)u^{-\alpha}$, which is integrable on $(0,\infty)$
($g'(u)\sim u^3/6$ at $0$, so $g'(u)u^{-\alpha}\sim u^{3-\alpha}$ with $\alpha<4$;
$g'(u)\le u$ at $\infty$, so $g'(u)u^{-\alpha}\le u^{1-\alpha}$ with $\alpha>2$).
Dominated convergence yields $J(t)=c_0\kappa_\alpha|t|^\alpha(1+o(1))$ with
$\kappa_\alpha=\int_0^\infty(u-\sin u)u^{-\alpha}du>0$ (the integrand is nonnegative
and not identically zero).
\end{proof}

\begin{lemma}[cf of $Z_n$ on a fixed annulus]\label{lem:cfZn}
Under \cref{ass:RV}, uniformly for $1\le|t|\le2$,
\[
\varphi_{Z_n}(t)=e^{-t^2/2}\Big(1+\kappa^\ast\,|t|^\alpha\,n^{1-\alpha/2}\big(1+o(1)\big)\Big),
\qquad \kappa^\ast:=\frac{c_0\,\kappa_\alpha}{\sigma^\alpha}>0 .
\]
\end{lemma}

\begin{proof}
For $n$ large, $u=t/(\sigma\sqrt n)$ is small uniformly on the annulus and
$\varphi(u)>\tfrac12$, so $\log\varphi(u)$ is defined and, by \cref{lem:cf} and
$\log(1-z)=-z+O(z^2)$ with $z=O(u^2)$,
\[
n\log\varphi\!\Big(\frac{t}{\sigma\sqrt n}\Big)
=-\frac{t^2}{2}
+\kappa^\ast|t|^\alpha n^{1-\alpha/2}\big(1+o(1)\big)
+O(n^{-1}),
\]
uniformly on $1\le|t|\le2$ (the $o(1)$ of \cref{lem:cf} depends only on the magnitude
of the argument, which tends to $0$ uniformly there). Since $\alpha<3$ implies
$1-\alpha/2>-\tfrac12>-1$, the $O(n^{-1})$ is $o(n^{1-\alpha/2})$. Exponentiating and
using $e^w=1+w(1+o(1))$ for $w\to0$ gives the claim; $\varphi_{Z_n}$ is real by
symmetry.
\end{proof}

\begin{theorem}[Impossibility]\label{thm:lower}
Assume \cref{ass:RV}. Let $w:\R\to(0,1]$ be Borel with $1/w$ of polynomial growth,
i.e.\ $w(x)\ge C_w^{-1}(1+|x|)^{-m}$ for some $m<\infty$. Then
\[
\liminf_{n\to\infty}\ n^{(\alpha-2)/2}\,
\sup_{t\in\R}\,w(t)\,\big|F_n(t)-\Phi(t)\big|\ \ge\ c(w,\alpha,c_0,\sigma)\ >\ 0 .
\]
In particular, this applies to $w=w_q=(1+\h)^{-q}$ for every $q\ge0$ and every $\h$
satisfying \cref{ass:h}.
\end{theorem}

\begin{proof}
Write $d_w:=\sup_t w(t)|F_n(t)-\Phi(t)|$, so that
\begin{equation}\label{eq:pointwise}
|F_n(x)-\Phi(x)|\ \le\ C_w\,(1+|x|)^m\, d_w\qquad(x\in\R).
\end{equation}
Choose $\chi\in C_c^\infty((1,2))$, $\chi\ge0$, $\chi\not\equiv0$, set
$\widehat\psi(t):=\chi(|t|)$ (even, real, smooth, compactly supported in
$\{1\le|t|\le2\}$) and
$\psi(x):=\frac{1}{2\pi}\int\widehat\psi(t)e^{itx}\,dt$, a real, even Schwartz
function. By Fubini (the $t$-integral is over a compact set),
\begin{equation}\label{eq:parseval}
\int_\R \psi\,d(F_n-\Phi)
=\frac{1}{2\pi}\int_\R \widehat\psi(t)\,\big(\varphi_{Z_n}(t)-e^{-t^2/2}\big)\,dt .
\end{equation}
By \cref{lem:cfZn}, the right-hand side equals
\[
\frac{\kappa^\ast}{2\pi}\,n^{1-\alpha/2}\,
\int_{1\le|t|\le2}\widehat\psi(t)\,|t|^\alpha e^{-t^2/2}\,dt\ \big(1+o(1)\big)
=:\frac{\kappa^\ast I_\psi}{2\pi}\,n^{1-\alpha/2}\,(1+o(1)),
\qquad I_\psi>0,
\]
using that $\widehat\psi\ge0$ and that the $o(1)$ in \cref{lem:cfZn} is uniform on the
support. For the left-hand side, integrate by parts: $\psi$ is Schwartz, $F_n-\Phi$ is
bounded and tends to $0$ at $\pm\infty$, so the boundary terms vanish and
\[
\Big|\int_\R \psi\,d(F_n-\Phi)\Big|
=\Big|\int_\R \psi'(x)\,\big(F_n(x)-\Phi(x)\big)\,dx\Big|
\ \le\ C_w\,d_w \int_\R |\psi'(x)|\,(1+|x|)^m\,dx\ =:\ C_1\,d_w,
\]
where $C_1<\infty$ because $\psi'$ is Schwartz, by \eqref{eq:pointwise}. Combining,
\[
d_w\ \ge\ \frac{\kappa^\ast I_\psi}{2\pi C_1}\; n^{1-\alpha/2}\,(1+o(1))
= c\, n^{-(\alpha-2)/2}\,(1+o(1)) . \qedhere
\]
\end{proof}

\begin{corollary}[Exact rate; the weight is rate-neutral]\label{cor:exact}
Under \cref{ass:RV,ass:h}, for every $q>0$ and every $\delta<\alpha-2$,
\[
c\,n^{-(\alpha-2)/2}\ \le\ d_{K,\h,q}\big(\mathcal L(Z_n),\Phi\big)
\ \le\ \sup_t|F_n-\Phi|\ \le\ C_\delta\, n^{-\delta/2},
\]
so both the weighted and the uniform Kolmogorov distances are
$n^{-(\alpha-2)/2+o(1)}$. In particular:
\begin{itemize}
\item[(i)] the rate $O(n^{-1/2})$ is unattainable for $\alpha\in(2,3)$, for every
choice of $(\h,q)$ and, by \cref{thm:lower}, for every weight with polynomially
bounded inverse;
\item[(ii)] for Student-$t$ with $\nu=2.5$, both metrics converge at
$n^{-1/4+o(1)}$ --- the rate the January version attributed to the uniform metric
alone.
\end{itemize}
\end{corollary}

\begin{remark}[Asymmetric tails]\label{rem:asym}
For non-symmetric regularly varying tails of index $\alpha\in(2,3)$ (e.g.\ one-sided
Pareto), the expansion of \cref{lem:cf} holds with a complex constant multiplying
$|t|^\alpha$ (plus a term $i\gamma t|t|^{\alpha-1}$), still nonzero; see
\cite[Ch.~2]{IbragimovLinnik} and \cite{ChristophWolf}. The proof of \cref{thm:lower}
goes through with $\widehat\psi$ testing the real or imaginary part, and the conclusion
is unchanged. We restricted to the symmetric case only to keep the argument fully
self-contained.
\end{remark}

\begin{remark}[What the weight \emph{does} do]\label{rem:constants}
By the non-uniform bound in \cref{thm:upper}, the contribution of the region
$\{|t|>T\}$ to $d_{K,\h,q}$ is $O(n^{-\delta/2}T^{-(2+\delta+q)})$: the weight
suppresses the tail region's share of the statistic. This is a statement about
\emph{constants} and about the finite-sample variance of the empirical statistic (a
single extreme observation moves the weighted statistic less than the uniform one),
not about the convergence rate. It is the honest version of the ``noise barrier''
intuition, and it is why the weighted statistic can still be a reasonable diagnostic
(\cref{sec:practice}).
\end{remark}

\section{What actually restores the rate}\label{sec:whatworks}

The lower bound identifies the obstruction exactly: the leading correction
$F_n-\Phi\approx \kappa^\ast n^{1-\alpha/2}\,G_\alpha$ is a fixed profile at bounded
argument. Two operations remove it; reweighting the metric is not one of them.

\begin{remark}[Correcting the target]\label{rem:target}
Define $G_\alpha$ by Fourier inversion of
$t\mapsto e^{-t^2/2}|t|^\alpha/(it)$, so that, under \cref{ass:RV},
$F_n=\Phi+\kappa^\ast n^{1-\alpha/2}G_\alpha+o(n^{1-\alpha/2})$ uniformly on compacts.
Comparing $\mathcal L(Z_n)$ with the corrected target
$\Phi_{n,\alpha}:=\Phi+\kappa^\ast n^{1-\alpha/2}G_\alpha$ removes the leading
obstruction \emph{by construction}. Quantifying the residual requires second-order
regular variation of the tail and is the subject of the Edgeworth-type theory for
distributions attracted to the normal law with heavy tails; see
\cite{ChristophWolf} and \cite[Ch.~2]{IbragimovLinnik}. We make no new quantitative
claim here.
\end{remark}

\begin{remark}[Trimming the data, not the metric]\label{rem:trim}
The truncation step inside the January version's own proof --- summing
$(X_i-\mu)\1\{\h(X_i)\le R_n\}$ --- is a winsorization of the \emph{data}, and it is
the part of the scheme that genuinely regains Gaussian behaviour: the truncated
variables have all moments, and the Berry--Esseen bound for them is honest. The scheme
failed only when it tried to transfer this back to the \emph{untrimmed} statistic
through the metric weight. The correct deployment is therefore to change the statistic:
trimmed and winsorized sums are asymptotically normal under regularly varying tails
with $\alpha>2$, with a well-developed distributional theory
\cite{CsorgoHaeuslerMason}, and validation can be performed on the trimmed statistic at
the classical rate, with the removed tail mass monitored by a separate tail diagnostic
--- precisely the hybrid structure of \cref{sec:practice}. The January version's
Remark~2.2, which dismissed winsorization in favour of smooth metric weighting, had the
moral inverted.
\end{remark}

\section{Numerical protocol: why the January figures were noise}\label{sec:numerics}

The January version estimated $d_{K,\h,q}(\mathcal L(Z_n),\Phi)$ by an empirical cdf
over $B=160$ Monte Carlo batches. For a sup-type statistic, the estimation error of an
empirical cdf based on $B$ samples is of order $B^{-1/2}$ (Kolmogorov:
$\E\sup_t|\widehat F_B-F|\approx 0.87/\sqrt B\approx 0.07$ for $B=160$), independently
of $n$. All quantities reported in Figures~1--2 of that version lie in the range
$0.03$--$0.12$ and are non-monotone in $n$: they are dominated by the Monte Carlo
floor, and the slopes stated in the captions (including the claim that $n=100$ matches
$n=10{,}000$) are not supported by the experiment. Those figures are withdrawn.

A correct protocol requires no Monte Carlo at all, and we implement it. $F_n$ is
computed by Gil--Pelaez inversion of $\varphi(t/(\sigma\sqrt n))^n$, with the
characteristic function evaluated through the cancellation-free decomposition
$1-\varphi(u)=\sigma^2u^2/2-J(u)$, $J(u)=\int g(ux)\,dF(x)$,
$g(y)=y^2/2-1+\cos y\ge0$ (adaptive quadrature of a nonnegative, non-oscillatory
integrand), and the exponent assembled as
$n\log\varphi(u)+t^2/2=nJ(u)+n\big(\log(1-\varepsilon)+\varepsilon\big)$,
$\varepsilon=1-\varphi(u)$, so that no precision is lost at large $n$. The pipeline is
cross-validated three ways: against an independent direct quadrature of $1-\varphi$
(relative agreement $10^{-8}$); against the prediction
$J(u)\sim c_0\kappa_\alpha u^\alpha$ of \cref{lem:cf} (the ratio $J(u)/u^\alpha$ is
constant to six digits over two decades of $u$); and against a $4\times10^5$-sample
Monte Carlo estimate of $F_n-\Phi$ at $n=256$ (agreement within one standard error at
each tested point).

\begin{figure}[H]
\centering
\includegraphics[width=0.74\textwidth]{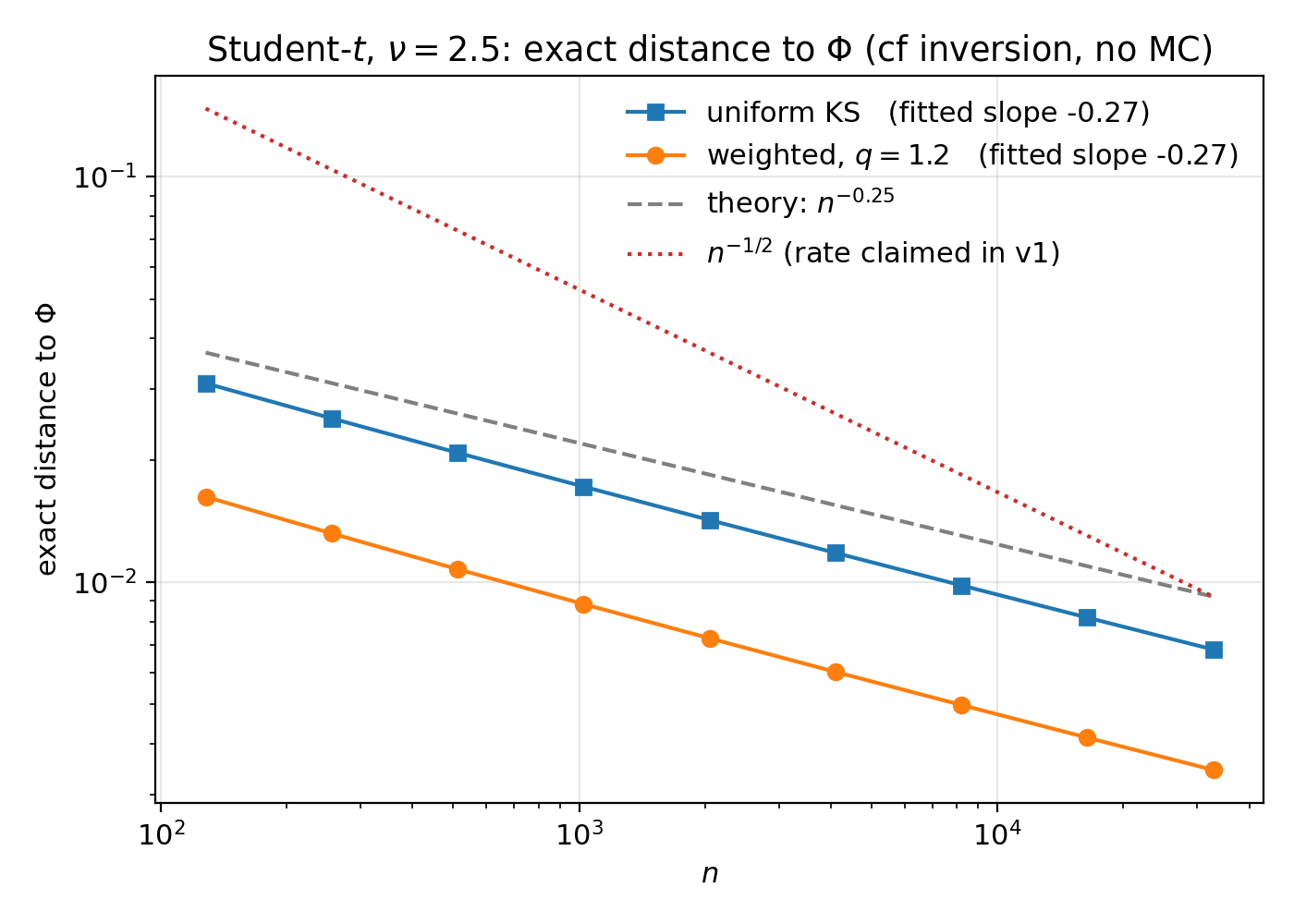}
\caption{Student-$t$, $\nu=2.5$: exact distances to $\Phi$, computed by deterministic
cf inversion. The uniform (squares) and weighted ($q=1.2$, circles) Kolmogorov metrics
decay at the same rate; fitted slopes $-0.27$ for both, theory $-(\nu-2)/2=-0.25$
(dashed). The dotted line is the $n^{-1/2}$ rate claimed in the January version.}
\label{fig:rate-t}
\end{figure}

\begin{figure}[H]
\centering
\includegraphics[width=0.74\textwidth]{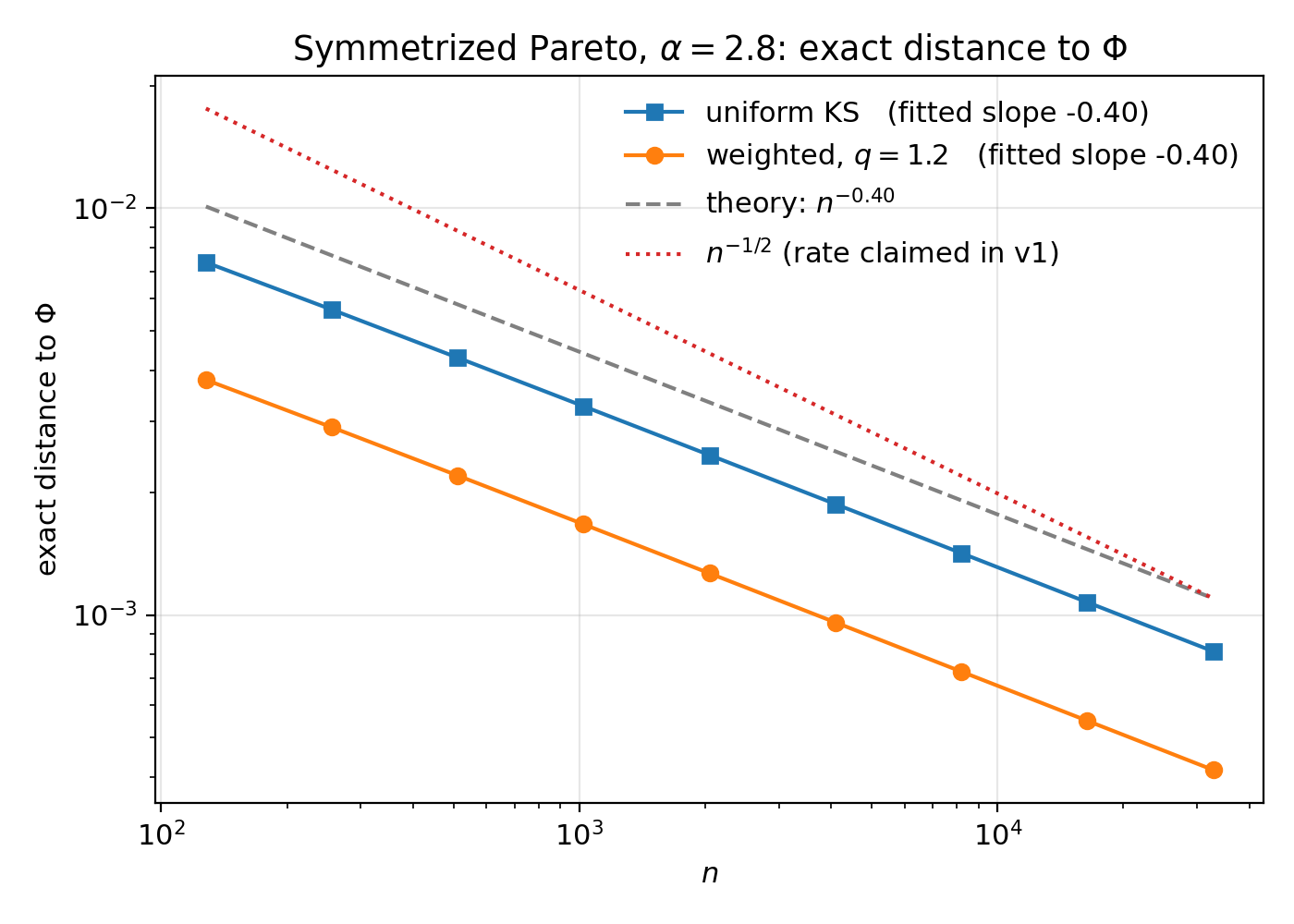}
\caption{Symmetrized Pareto, $\alpha=2.8$: fitted slopes $-0.400$ (uniform) and
$-0.401$ (weighted) against the theoretical $-(\alpha-2)/2=-0.40$. The weight shifts
the constant, not the rate.}
\label{fig:rate-p}
\end{figure}

\begin{figure}[H]
\centering
\includegraphics[width=0.74\textwidth]{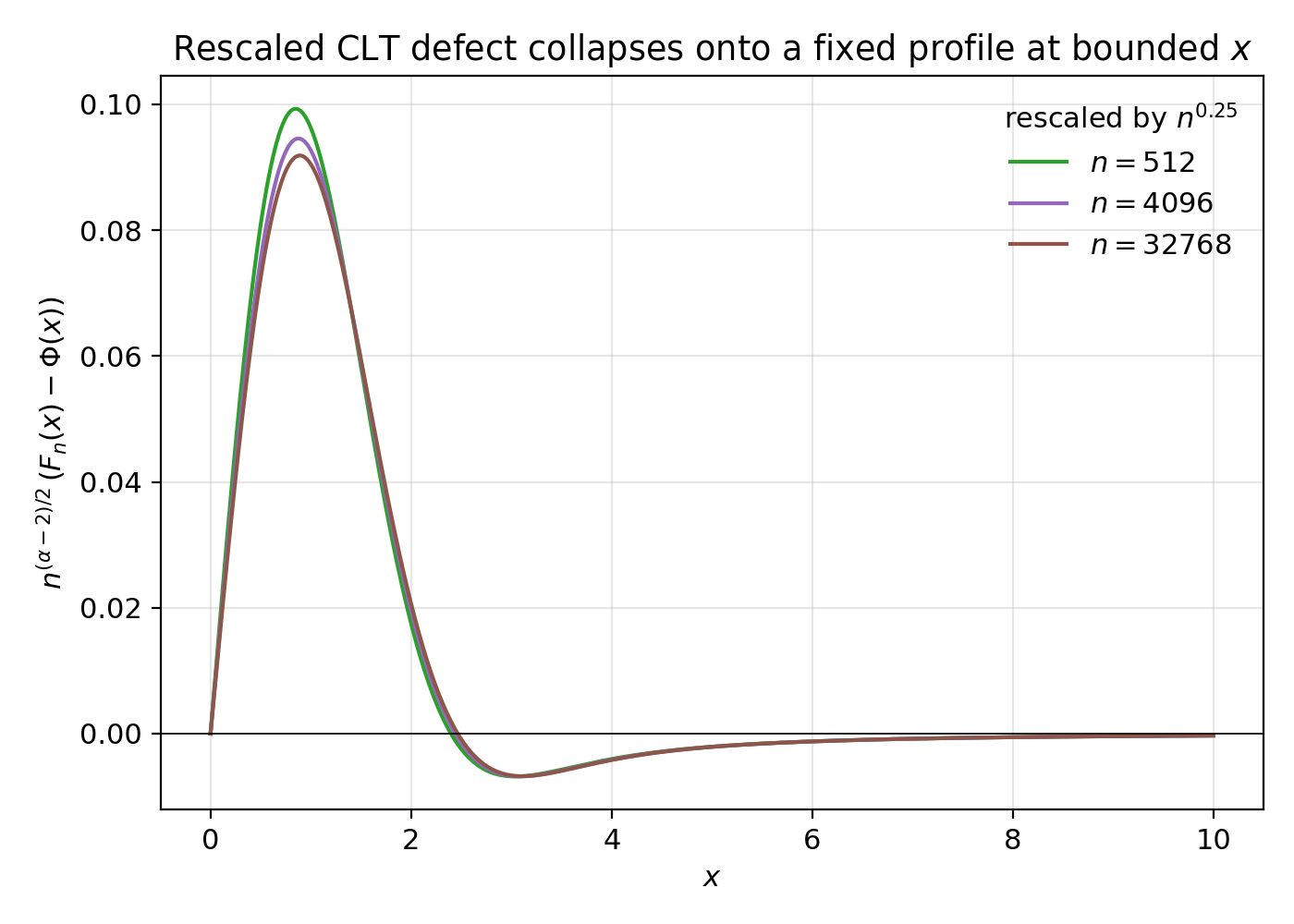}
\caption{Rescaled defect $n^{(\alpha-2)/2}(F_n-\Phi)(x)$ for $n=512,4096,32768$
(Student-$t$, $\nu=2.5$): the curves collapse onto a fixed profile located at bounded
argument, where every weight with polynomially bounded inverse is $\asymp1$ --- the
geometric content of \cref{thm:lower}.}
\label{fig:collapse}
\end{figure}

\begin{figure}[H]
\centering
\includegraphics[width=0.74\textwidth]{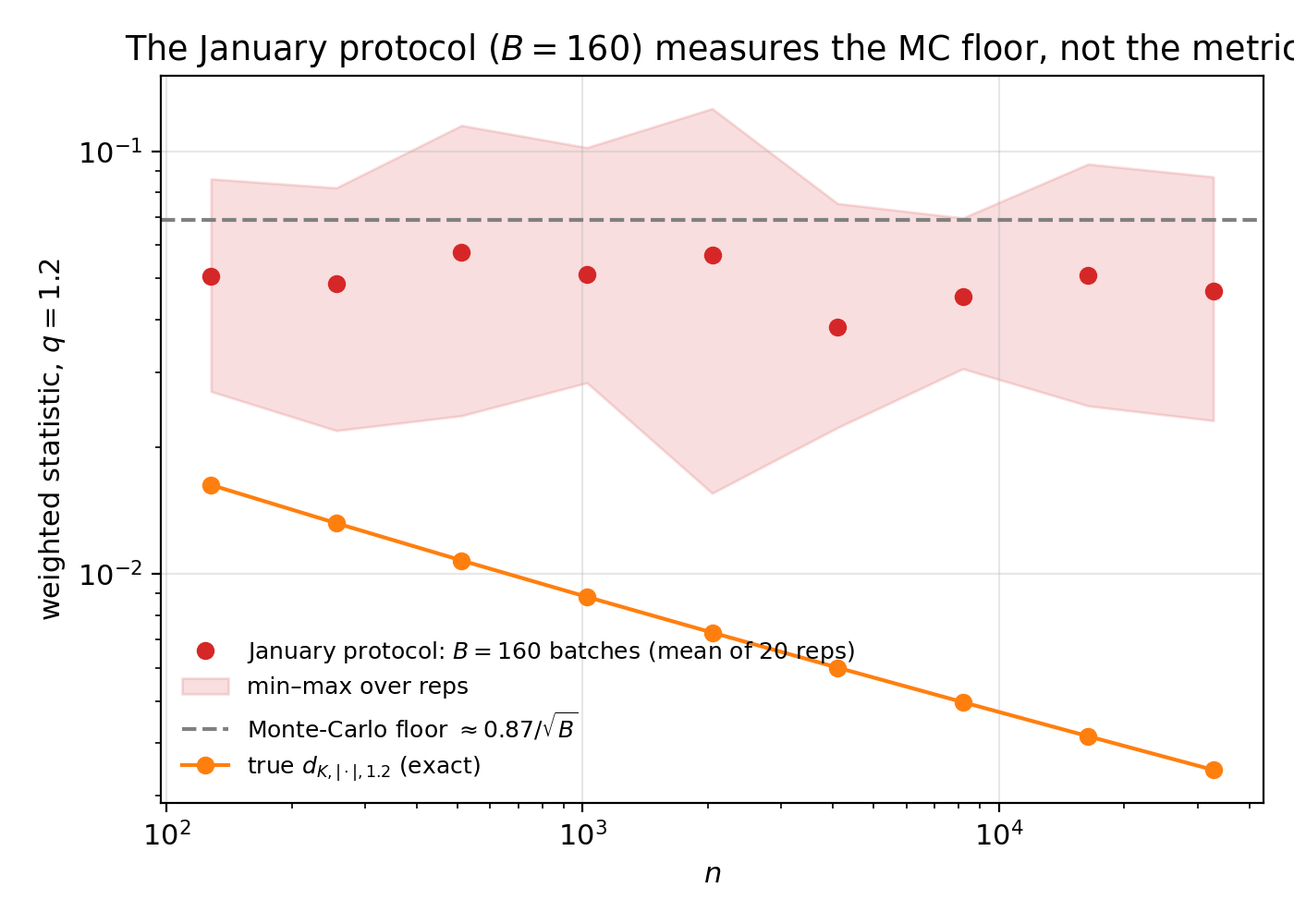}
\caption{The January protocol reproduced: the $B=160$ batch statistic (dots, $20$
repetitions, min--max band) is flat in $n$ at the Monte Carlo floor (dashed), an order
of magnitude above the exact weighted metric (solid). The non-monotone wiggles of the
January figures were sampling noise, not a convergence rate.}
\label{fig:floor}
\end{figure}

\Cref{fig:rate-t,fig:rate-p} display the exact $d_K^{\mathrm{unif}}$ and
$d_{K,|\cdot|,1.2}$ for $n=2^7,\dots,2^{15}$. Fitted log-log slopes: Student-$t$
($\nu=2.5$): $-0.27$ for \emph{both} metrics, against the theoretical
$-(\nu-2)/2=-0.25$ (the residual curvature is the slowly decaying next-order
correction); symmetrized Pareto ($\alpha=2.8$): $-0.400$ and $-0.401$, against the
theoretical $-0.40$. In both cases the two curves are parallel: the weight divides the
constant by roughly two and leaves the rate untouched, exactly as \cref{cor:exact}
states, and both are visibly slower than the $n^{-1/2}$ guide claimed in the January
version. \Cref{fig:collapse} shows the mechanism of \cref{thm:lower} directly: the
rescaled differences $n^{(\alpha-2)/2}(F_n-\Phi)$ for $n=512,\,4096,\,32768$ collapse
onto a single fixed profile, peaked at $x\approx0.9$ and supported on $|x|\lesssim4$
--- the region where every admissible weight is of order one. Finally,
\cref{fig:floor} re-runs the January protocol ($B=160$ batches, $20$ repetitions)
against the exact curve: the empirical statistic is flat in $n$ at the Monte Carlo
floor $\approx0.87/\sqrt B\approx0.07$ (means $0.04$--$0.06$ for the weighted
version), an order of magnitude above the true metric, and carries no information
about any rate. Resolving the $n^{-1/2}$-vs-$n^{-(\alpha-2)/2}$ distinction at
$n=32{,}768$ by Monte Carlo would require $B\gtrsim10^5$ batches. The reproduction
script is distributed with this note.

\section{What survives for practice}\label{sec:practice}

Stripped of the rate claim, the weighted statistic remains a legitimate object: any
functional of $(\widehat F_n,F_\theta)$ can serve as a test statistic once its null
distribution is calibrated. The operational material of the January version survives in
this reframed form, and we summarise it.

\emph{Calibration.} Critical values of $d_{K,\h,q}$ are model-dependent; compute them
by parametric bootstrap under $H_0$ (simulate from the fitted $F_\theta$, take empirical
quantiles of the bootstrap statistics). \emph{Hybrid acceptance.} Because the weight
deliberately reduces tail sensitivity (\cref{rem:constants}), the core statistic must
be paired with a dedicated tail safeguard (VaR exception counts, Kupiec/Christoffersen,
or an ES backtest): accept iff both pass. \emph{Anti-gaming.} Pre-specify $q$ (or a
grid $Q$, accepting only if all $q\in Q$ pass) and separate calibration and validation
windows. \emph{Motivation, restated honestly.} The gain from the weight is
finite-sample robustness of the diagnostic --- a single extreme observation perturbs
the weighted statistic by $O(w_q(x_{\max})/n)$ rather than $O(1/n)$ uniformly --- not a
faster rate of convergence, which \cref{thm:lower} rules out.

\section{Conclusion}

Under sub-cubic moments the Berry--Esseen defect is not tail noise contaminating an
otherwise Gaussian centre: it \emph{is} a deficiency of the centre, of size
$n^{-(\alpha-2)/2}$, produced by the $|t|^\alpha$ term of the characteristic function.
No reweighting of the Kolmogorov distance in the argument can remove it
(\cref{thm:lower}); the exact rate of the weighted metric coincides with the uniform
one (\cref{cor:exact}). The rate can be restored only by correcting the Gaussian
target or by trimming the statistic (\cref{sec:whatworks}); the weighted metric itself
is at best a variance-reduced diagnostic to be calibrated by bootstrap
(\cref{sec:practice}). The January 2026 version of this note, which claimed the
opposite, is withdrawn; \cref{app:corrections} lists the corrections item by item.

\appendix

\section{Corrections with respect to the January 2026 version}\label{app:corrections}

\begin{enumerate}
\item \textbf{Theorem 3.9, Corollary 4.4 (weighted BE at $n^{-1/2}$).} False. The
conclusion is contradicted by \cref{thm:lower} for every regularly varying tail of
index $\alpha\in(2,3)$, i.e.\ on the paper's own headline examples (Student-$t$
$\nu=2.5$, Pareto $\alpha=2.8$). Irreparable within the class of weighted sup-metrics:
the failure is of the statement, not of the proof.
\item \textbf{Proposition 4.2 (choice of $(\beta,q)$).} Arithmetically wrong: the
condition $\beta(1-\delta)\le\frac12$ makes the core term $O(1)$, not $O(n^{-1/2})$;
the region compatible with $n^{-1/2}$ forces $\delta=1$ (\cref{prop:empty}). The
recommended $\beta^\star=1/(2\eta)$ makes the bound divergent when $\eta<1-\delta$
(\cref{rem:arith}). Theorem 3.9 moreover imposed no constraint on the core term at
all.
\item \textbf{Lemma 4.1 / Theorem 3.5 (core/tail scheme).} The truncation remainder
omits the factor $n$ in $n\,P(\h(X)>R)$ and the mean-shift term
$\sqrt n\,|\E[(X-\mu)\1\{\h>R\}]|$; the constants cannot be absolute as stated
(\cref{rem:missing-n}). Superseded by \cref{thm:upper}, which is stronger and does not
use truncation.
\item \textbf{Theorem 6.2 (multivariate version).} Retracted for the same two reasons;
the univariate lower bound already refutes the rate claim by embedding.
\item \textbf{Remark 3.10 and the abstract/introduction.} The claims that the
$n^{-1/2}$ rate is ``restored'' for Pareto and Student-$t$, and the condition
``$q\ge\eta$'', are withdrawn.
\item \textbf{Section 7 (numerical validation) and Figures 1, 2, 4.} Withdrawn. With
$B=160$ batches the Monte Carlo floor ($\approx0.07$ for sup-statistics) exceeds every
effect the figures purport to measure; the stated slopes and the ``$n=100$ vs
$n=10{,}000$'' equivalence are artifacts (\cref{sec:numerics}). The correct prediction
is that both metrics decay at $n^{-(\alpha-2)/2}$.
\item \textbf{Remark 2.2 (winsorization).} Moral inverted: truncating the data is
precisely the operation that restores Gaussian rates; reweighting the metric is the
operation that achieves nothing (\cref{rem:trim}).
\item \textbf{Kept.} Proposition 2.4 (metric property), Lemma 3.3 / Proposition 5.1
(weight and metric equivalence under coarse change of exhaustion), Proposition 3.8
(tail remainder under regular variation), and the operational material of Section 8,
reframed without the rate claim in \cref{sec:practice}. \Cref{ass:h} now states local
boundedness explicitly, which the January version used implicitly.
\end{enumerate}

\end{document}